\documentclass[12pt]{amsart}
\usepackage{amssymb}
\begin{document}

\title{Counting Singular Plane Curves Via Hilbert Schemes}
\author{Heather Russell}

\maketitle
\tableofcontents
\date{\today}
\def\Q{{\mathbb Q}}
\def\A{{\mathbb A}}
\def\P{{\mathbb P}}
\def\C{{\mathbb C}}
\def\G{{\mathbb G}}
\def\F{{\mathbb F}}
\def\Z{{\mathbb Z}}
\def\Sp{{\textup{Spec}}}
\def\O{{\mathcal O}}
\def\H{{\textup{Hilb}}}
\def\N{{\mathbb N}} 
\def\I{{\mathcal I}}
\def\m{{\mathfrak m}}
\def\Se{{\mathcal S}}
\def\aut{{\rm Aut}}
\def\L{{\mathcal L}}

\newtheorem{example}{Example}[section]
\newtheorem{theorem}{Theorem}[section]   
\newtheorem{corollary}{Corollary}[section]
\newtheorem{question}{Question}
\newtheorem{lemma}{Lemma}[section]
\newtheorem{proposition}{Proposition}[section]
\section{Introduction}

Consider the following question.  Given a suitable linear series $\L$ 
on a smooth surface $S$, how many 
curves in $\L$ have a given analytic or topological type of 
singularity?   
By ``suitable'' linear series with respect to a type of singularity, we
mean that there are finitely many curves with the singularity in 
the linear series and their codimension is maximal.  
Our approach to this question is to express the answer as a Chern number
of a vector bundle over a compactification of 
a space linearizing the condition of having
the singularity.  For example, the condition of having a cusp along a
given tangent direction at a given point is linear in the sense that
curves in a linear series spanned by two curves with the 
condition also have the condition.  Thus, 
the projectivized tangent bundle 
$\P T(S)$ linearizes the condition of having a cusp in $S$.  Note 
that the closure of the condition of 
having a cusp in along a given direction is
the condition of containing a particular subscheme of $S$ isomorphic
to $\Sp (R/(x^2,xy^2,y^3))$ where $R$ is the ring $K[[x,y]]$ and $K$
is the field of definition of $S$.  Generalizing this, the spaces we will
use to linearize our conditions will be of the form
$$U(I)=\{a \in \H ^d (S): a \cong {\rm Spec}(R/I)\}$$ for 
ideals $I$ of some finite colength $d$ in $R$.        
Section 3 is devoted to this correspondence between ideals and types 
of singularities and other conditions on curves. 
The space $U(I)$ has a natural compactification $C(I)$, its
closure in $\H ^d (S)$.  Letting $L$ be the line bundle corresponding
to the divisor of a section of $\L$,  $C(I)$ admits the vector bundle 
$V_{L}(I)$ defined in section 4 with the property that sections 
of $\L$ give sections of the $V_{L}(I)$ vanishing exactly over those
points corresponding to the data of the singularity in question of the
curve corresponding to the section.  Thus the number of curves having
the given types of singularity is the number of places a set of
sections of $V_{L}(I)$ coming from a basis of sections of $\L$ become
dependent.  This is the Chern number of $\L$.  In order to express 
this Chern number in terms of the divisor $D$ of a section of $L$ and
the Chern classes of the tangent bundle of $S$, one would like to know
the Chern polynomial of $V_{L}(I)$ pulled back to some space with a
known Chow ring, so that the relations in the Chow ring can be used to
find the degree of the relevant Chern class.  In the case that 
$I$ can be constructed from the ideals $(x,y^3)$ and $(x,y)$ or from 
$(x,y^2)$ and $(x^2,y)$ by taking sums, products, and images under the 
Frobenius morphism if we are working over positive characteristic we can 
both find the Chow ring of a space dominating $C(I)$ 
and the Chern classes of $V_L(I)$ in terms of this Chow ring.  However, 
if $I$ is constructed from the ideals $(x,y^4)$ and $(x,y)$ or $(x,y^3)$ 
and $(x^2,y)$, we can find the Chow ring of a space dominating $C(I)$, 
but do not have an algorithm for finding the Chern classes of $V_L(I)$. 
However, by slightly ad-hoc means one can sometimes or possibly always
find these Chern classes.  In the last section we give examples of both
enumerative problems solved solely by previous results and enumerative
problems solved by a mixture of previous results and ad-hoc means.  
                  
\vspace{.1 in}

\noindent {\bf Acknowledgements:}  I would like to thank Karen
Chandler, Joe Harris, Tony Iarrobino, Steve Kleiman, Suresh Nayak, 
Dipendra Prasad, 
Mike Roth, Jason Starr, Ravi Vakil, Joachim Yameogo, 
and many others for their generous
help.  

\section{Preliminaries}  The foundation of this paper has been built
up in \cite{R1} and \cite{R2}.  We recall some definitions and
notation.  

We will use the following short hand for denoting monomial ideals.  
Given a sequence of positive integers 
$s_1, \dots, s_r$, we will let $I(n_1, \dots, n_r)$ denote the ideal 
$(x^r, x^{r-1}y^{n_1}, x^{r-2}y^{n_1+n^2}, \dots )$. 
Moreover, given a sequence of monomial ideal, $I_1, \dots , I_r$ we
let 
$$U(I_1, \dots, I_r) = 
\{ (a_1, \dots ,a_r) \in U(I_1) \times \dots \times U(I_r):  
\exists p \in X ~{\rm and} 
$$ $$ \varphi : R \tilde{\rightarrow} {\hat \O}_{X,p} ~{\rm with}~ 
\varphi (I_1, \dots , I_r) = (a_1, \dots ,a_r)\}$$ 
and $C(I_1, \dots, I_r)$ be its closure in the appropriate product of Hilbert 
schemes.  We will say that 
$C(I_1, \dots,I_r)$ is an {\it alignment correspondences} with 
{\it interior} 
$U(I_1, \dots , I_r)$.  Moreover, we will say that the measuring
sequence of  $I_1, \dots , I_r$ is $A_1$, $A_2$ where $A_1$
(respectively $A_2$) is the
ideal generated by images of $x$ (respectively $y$) under 
automorphisms of $R$ fixing $y$ (respectively $x$) and sending the
$I_j$'s to themselves.  
In accordance with \cite{R1}, although not \cite{R2}, we let 
$G(I_1, \dots , I_r)$ be the group of automorphisms of $R$ sending the
$I_j$'s to themselves.

\section {Conditions on Curves}

\noindent{\bf Definition}: Say that a curve $C$ has the {\it condition} 
corresponding to an ideal $I$ if $C$ contains a subscheme 
corresponding to a point of $C(I)$.  

\vspace{.1 in}

We will use the following  
two lemmas to identify conditions corresponding to ideals.
\vspace{.1 in}

\vspace {.1 in}

\noindent{\bf Definition}: Say that a curve $C$ is {\it generic} with 
the condition corresponding to an ideal $I$ if $I$ imposes independent 
conditions on the complete linear series containing $C$ and $C$ is generic 
among curves in this linear series with this condition.    

\begin{lemma}\label{transform}  Given an ideal $I$ of finite colength in 
$R$, the condition on the proper transform of 
a generic curve with the condition corresponding to $I$ with respect to the 
blow up of $S$ at the singular point of the curve is the condition 
corresponding to the quadratic transform of $I$ as defined in \cite{Z}.    
In particular, given a sequence of positive integers $n_1, \dots ,n_r$, 
the condition corresponding to $I(n_1, \dots ,n_r)$ is the closure of the condition of having an $r^{\rm th}$ order point such that 
if $y$ is a local coordinate for the exceptional divisor of 
the proper transform of a generic 
curve $C$ with the condition, then the condition on the proper 
transform corresponds to the ideal 
$I(n_1 -1, \dots , n_r -1)$.
\end{lemma}

\noindent{\bf Proof:} The lemma can be verified by direct computation in 
coordinate patches.  $\square$

\begin{lemma}\label{integralclosure}  The condition 
corresponding to the integral closure of an ideal 
$I$ is the closure of the condition of having the topological type of 
singularity of a generic curve with the condition corresponding to $I$. 
\end{lemma}

\noindent{\bf Proof:} By \cite{Z} the quadratic transform of 
integral closure of $I$ is the integral closure of the quadratic transform 
of $I$.  Thus the lemma follows from Lemma~\ref{transform} and induction 
on the colength of $I$. $\square$

\vspace{.1 in}
\begin{theorem}\label{degenerationideal} Let $I$ be an ideal with 
measuring sequence at most $(x,y^3),(x,y)$ 
(respectively $(x,y^2), (x^2,y)$).  Then the boundary of $C(I)$ 
is equal to the space $C(J)$ where 
$$J = lim _{t \to \infty}g(t) (I)$$
and $g(t)$ is the automorphism of $R$ sending $x$ to $x+ty^2$ and fixing 
$y$ (respectively fixing $x$ and sending $y$ to $y+tx$.)
\end{theorem}

\noindent{\bf Proof:}  The automorphisms of the form $g(t)$ form a set of 
coset representatives of $G((x,y^2)/G(I, (x,y^2))$.  Therefore the 
boundary of the fiber of 
$C(I, (x, y^2))$ over $C((x,y^2))$ is the fiber of 
$C(J, (x,y^2))$ over 
$C((x,y^2))$ with respect to the projection map.  Hence projecting 
the space $C(I,(x,y^2))$ to $C(I)$, the proposition follows.     
$\square$   

\vspace{.1 in}

\noindent{\bf Definition:} Given an ideal $I$ as in
Theorem~\ref{degenerationideal}, we will say that $J$ is the
degeneration ideal of $I$.  

\vspace{.1 in}

\noindent{\bf Definition:} We will call the  
codimension of curves with the condition given by an ideal $I$ 
in a linear series on which $I$ imposes independent conditions 
is {\it codimension} of the condition, denoted ${\rm cod}(I)$.  

\begin{lemma}\label{equality}The codimension 
of the condition corresponding to $I$ satisfies 
$${\rm cod}(I) + {\rm dim}(C(I))= {\rm col}(I)  +\epsilon(I)$$
where $\epsilon(I)$ is the dimension of the locus of subschemes in 
$C(I)$ contained in a generic curve with the condition given by $I$.  
\end{lemma}  

\noindent{\bf Proof:} Let $\L$ be a linear series for which $I$ imposes 
independent conditions.  Consider the incidence correspondence  
$$\{(C, \alpha )\in \L \times C(I):  \alpha \subset C\}.$$  
Equating sums of the dimension of the base and fiber with respect to the 
two projection maps, the lemma follows. $\square$    

\vspace{.1 in}

\begin{lemma}\label{versal} The codimension of the condition of 
having an analytic type of singularity 
or a degeneration is the dimension of the  
versal deformation space of that singularity.          
\end{lemma}

\noindent{\bf Proof:} The versal deformation space to a 
singularity can be naturally identified with the normal space to 
the tangent space of the locus of curves in a linear series such that there 
is an ideal imposing independent condition on the linear series and 
generic members have the given type of analytic singularity.  $\square$
 
\begin{example} \end{example} 
Let $n_1, \dots , n_r$ be a sequence of increasing positive integers.  
Then from  
Lemma~\ref{transform} one can see that the 
condition of the condition given by the ideal $I(n_1, \dots , n_r)$ 
is the closure of the topological condition with 
Enriques diagram a succession of $r$ free vertices 
of decreasing weights such that there are exactly $n_{i}$ vertices of 
weight at most $r+1-i$.  

\begin{example} \end{example}  A generic curve with the condition 
given by the ideal $(x^a,y^b)$ has the singularity of topological type 
$x^a+y^b$.  The integral closure of the ideal $(x^a,y^b)$ is the ideal 
$I(a,b)$ generated by monomials $x^cy^d$ with $3c + d \ge b$.  Hence 
this ideal gives the closure of the condition of having the singularity 
of topological type $x^a+y^b$ and has measuring sequence 
$(x,y^{[\frac {b}{a}]}), (x,y)$. 
  
\begin{example}\end{example} The condition given by the ideal  
$I(2,2,1,1)$ is strictly in between the closures of the 
conditions of having topological type and analytic 
type $x^4 +y^6$.  The topological condition is that of having a double
cusp.  If one blows up at the double cusp of a curve, 
the proper transform will be a tacnodal curve with both branches
tangent to the exceptional divisor.  If one then blows up at this
point of intersections one gets a quadruple point, with two branches
corresponding to the proper transform of the proper transform and two
branches corresponding to exceptional divisors.  The additional
condition corresponding to $I(2,2,1,1)$ is that the four tangent
directions of these four branches have cross-ratio $-1$.  Such curves 
are of analytic type $x^4+ax^2y^3+y^6$ for some constant $a$.  The
analytic type varies with $a^2$.  Unlike the analytic condition 
above given by a cross ratio, the analytic condition
corresponding to a particular choice of $a^2$ cannot be realized by
the configuration of points of intersection of components of the 
total transform of the curve after some number of blow-ups.  It would
be interesting to find a geometric way of visualizing 
such analytic conditions.  

\vspace{.1 in}

\begin{example} \end{example} Theorem~\ref{degenerationideal} 
can be used to glean some information about 
degenerations of singularities.  Given an ideal $I$ as in 
Theorem~\ref{degenerationideal}, the condition given by 
$I$ is that of containing either $\Sp (R/I)$ or $\Sp (R/J)$ 
where $J$ is the degeneration ideal of $I$.  

For example, degeneration ideal of $I(3,3,2)$ is the ideal 
$I(2,2,1,1)$.  This corresponds to 
the fact that the singularity of topological type $x^3 +y^8$ can
degenerate to the singularity of topological type $x^4 +y^6$.
However, only those curves with the additional analytic condition
described in the example above can occur.

\vspace{.1 in}

The sequence giving a condition is not in general unique.  For 
example, the condition corresponding to $1,n$ for any positive 
integer $n$ is that of being singular.  
Note that for $n>2$, a curve with a node will contain 
two schemes corresponding to points in $C((x^2, xy, y^2))$.  
The following proposition gives a criterion for when two ideals 
correspond to the same condition. 

\begin{proposition}  
Let $I_1$ and $I_2$ be ideals in $R$ with 
$I_1 \supsetneq I_2$.   If 
$${\rm cod}(I_1)= {\rm cod}(I_2), $$ 
then both ideals give the same condition.  
\end{proposition}

\noindent{\bf Proof}:  Let be $\L$ be the linear corresponding to a
sufficiently high tensor power of an ample line-bundle on $S$.  For 
$i \in \{1,2\}$, let $\Gamma_i$ be the incidence correspondence in $\L \times
C(I_i)$ as in Lemma~\ref{equality}.  Then each $\Gamma _i$ is a vector
bundle over $C(I_i)$ and hence irreducible.  Therefore the image of
the projection map $\pi_1: \Gamma_i \rightarrow \L$ is irreducible.
Since $\pi _1(\Gamma_1)$ contains $\pi _2(\Gamma_2)$ and these two
images are both of the same codimension in $\L$, they must be equal.
Suppose by way of contradiction that the $G(I_1)$ orbit of $I_2$ is
not equal to $I_1$.  Then there is an element $a \in I_1$ that is not
in the $G(I_1)$ orbit of $I_2$.  One can find a curve $C \in \L$ with
local equation at a point $p \in S$ equal to the image of $a$ under an
isomorphism from $R$ to $\hat \O_{S,p}$ up to
an element in the image of a high power of the maximal ideal in $R$.  
Then $C$ is in  
$\pi _1(\Gamma_1)$ but not $\pi _2(\Gamma_2)$.  It follows that 
$G(I_1)$ orbit of $I_2$ is $I_1$ and hence that $I_1$ and $I_2$
correspond to the same conditions. $\square$

\vspace{.1in}

The converse of the previous proposition does not hold in general
because one can have ideal corresponding to the same condition,
neither of which is contained in the other.  However, both will be
contained in the ideal of minimal colength corresponding to the condition. 

\section{Vector Bundles on Alignment Correspondences}

For each positive integer $d$ let $${\mathcal U} \subset \H^d(X)\times X$$ 
be the universal family over $\H^d(X).$ 

Let 
$$\pi:{\mathcal U}\rightarrow \H^d(X)$$
and 
$$\mu:{\mathcal U} \rightarrow X$$
be the projection maps.

\vspace{.1 in}

\noindent{\bf Definition:}  
Given a line bundle $L$ on $S$, let $V_L(I)$  
denote the restriction of $(\pi)_*\mu^*(L)$ to $C(I)$ or by 
abuse of notation its pullback to any
space mapping to $C(I)$.  

\vspace{.1 in}

The fiber in $V_L(I)$ over a point in $C(I)$
is the vector space of germs of sections 
of $L$ modulo those in the ideal corresponding to the point.  If $\L$
is a linear series such that $L$ is the line bundle corresponding to
the divisor of a section, then global sections of $\L$ 
give global sections of $V_L(I)$.  
These sections of $V_L(I)$ vanish exactly over the points parametrizing 
subschemes of the divisor of the corresponding section of $\L$.  Hence 
such a section vanishes if and only if the curve has the condition 
corresponding to $I$. 

By the following lemma, to find the Chern polynomial of $V_L(I)$, 
it is enough to find the 
Chern polynomial of $V(I)$.

\begin{lemma}\label{tensor}  
Let $A$ be a space admitting the vector bundle $V(I)$.  Then 
for any line bundle $L$ on $S$, we have an equality of Chern polynomials
$$c(V_L(I))=c(V(I)\otimes (L))$$ over $A$.  Here, abusing notation, 
we use $L$ to denote its pullback to $A$.
\end{lemma}

\noindent{\bf Proof}: Let $A^{\prime}$ be an extension $A$ over which
$V(I)$ has a filtration  
$$V(I)=V(I_1)\rightarrow \dots \rightarrow
V(I_r)\rightarrow 0$$
such that the successive kernels are line bundles.  Then 
$$V_L(I) \rightarrow V_L(I_1)\rightarrow \dots \rightarrow
V_L(I_r)\rightarrow 0$$
is a filtration of $V_L(I)$.  
From the isomorphisms
$${\rm Ker}(V_L(I_j) \rightarrow V_L(I_{j+1}))
\cong {\rm Ker}(V(I_j) \rightarrow V(I_{j+1}\otimes L))$$
and the Whitney product formula, we see that the lemma holds over 
$A^{\prime}$.  Thus it also holds over $A$.  
$\square$

\begin{proposition}\label{proj} Let $C(I_1, \dots, I_r)$ 
be an alignment correspondence 
such that there is a monomial ideal $I$ with 
$$I_2 \subset I \subset I_1$$ and 
$${\rm dim}(C(I_1, \dots, I_r, I))-{\rm dim}(C(I_1, \dots, I_r))= 
{\rm dim}(I_1/I_2)-1.$$ 
If $I$ has dimension $1$ (respectively codimension $1$) as a 
subspace of $I_1/I_2$, then the space $C(I_1, \dots, I_r, I)$ is the 
projectivization of the vector bundle $V(I_1/I_2)$ 
(respectively $V(I_1/I_2)^*$) over $C(I_1, \dots ,I_r)$.  

\end{proposition}

\noindent {\bf Proof:} If $I$ has dimension $1$ 
(respectively codimension $1$) as a 
subspace of $I_1/I_2$, then $C(I_1, \dots , I_r, I)$ 
has a natural embedding in $\P V(I_1/I_2)$ (respectively $\P V(I_1/I_2)^*$).  
Since both spaces are irreducible and of the 
same dimension, this embedding is an isomorphism.  $\square$  

\vspace{.1 in}

The Chow rings of the spaces $C((x,y^2),(x,y^3))$ and 
$C((x,y^2),(x^2,y))$ will be particularly useful for enumerative applications 
due to the fact that they are universal fiberwise $\aut (R)$-equivariant 
compactifications of the spaces $U((x,y^3))$ and $U((x,y^2),(x^2,y))$
respectively.  

\begin{lemma}\label{basering}\cite{Col}  
Let $I_i$ denote the ideal $(x,y^i)$.  Let $c_1$ and $c_2$ denote 
the first and second Chern classes of the cotangent bundle of $S$, 
respectively.  Let $h_2$ denote the hyperplane class of the projectivization 
of the cotangent bundle of $S$.  
Let $h_3$ and $h_3^{\prime}$ be the hyperplane class of 
the spaces $C(I_2,I_3)$ and $C(I_2,(x^2,y))$ 
as the projectivization of the bundles $V(I_2/I_1I_2)$
and $V(I_1/I_1^2)$ respectively.  Their Chow rings are given by

$$A(C(I_2,I_3)) = A(B)[h_3]/(h_3+2h_2+2c_1)(h_3-h_2)$$
and   
$$A(C(I_2,(x^2,y)) = 
A(B)[h_3^{\prime}]/((h_3^{\prime})^2+h_3^{\prime}c_1+c_2).$$
\end{lemma}  

\noindent{\bf Proof:}  The expression for the second Chow ring follows from 
the fact that $C((x,y^2)(x^2,y)$ is the projectivization of the
pullback of the cotangent 
bundle on $S$ over $B$.  

To find the Chow ring of $C(I_2,I_3)$ we will use the following 
two short exact sequences of 
bundles.       
\begin{equation}
0\rightarrow V(I_2/I_1^2) \rightarrow 
V(I_1/I_1^2)\rightarrow V(I_1/I_2)\rightarrow 0
\end{equation}
\begin{equation}
0\rightarrow V(I_1/I_2)^2 \rightarrow
V(I_2/I_1I_2)\rightarrow 
V(I_2/I_1^2) \rightarrow 0
\end{equation}
Since the vector bundle $V(I_2/I_1^2)$ 
is the tautological bundle over $B,$ 
$$c_1(V(I_2/I_1^2))=-h_2.$$  Applying the Whitney product formula 
to the two above sequences, we get 
$$c_1(V(I_2/I_1^2)) = c_1+h_2$$
and thus  
$$c(V(I_2/I_1I_2))= (1+2h_2+c_2)(1-h_2),$$ giving us the presentation of 
$A(C(I_2,I_3))$ as claimed. $\square$

\begin{proposition}\label{shift} 
Let let $V(I/J)$ be a bundle of rank one defined on $C((x,y^3))$ 
(respectively, $C((x,y^2), (y,x^2))$ where $I$ and $J$ are monomial
ideals such that there quotient is generated by    
$x^ay^b$ as a vector space.  Let $x^cy^d$ be 
the monomial generating the degeneration ideal of $I$ over the degeneration 
ideal of $J$.  
  
Then we have  
$$c_1(V(I/J))= 
-ah_2 + b(c_1+h_2) +(a-c)(h_2-h_3)$$
(respectively 
$$-ah_2 + b(c_1+h_2) + (a-c)(c_1+h_3^{\prime}+h_2)).$$ 

\end{proposition}

\noindent{\bf Proof}:  Let $I_j$ be the ideal $(x,y^j)$.  
Let $I^{\prime}$ and $J^{\prime}$ be 
the maximum monomial ideals 
with respect to inclusion with measuring sequence at most 
$I_3, I_1$ (respectively $I_2,(x^2,y)$) such that $I^{\prime}$ 
is generated over $J^{\prime}$ by $x^ay^b$ and the degeneration ideal of 
$I^{\prime}$ is generated by $x^cy^d$ over the degeneration ideal of 
$J^{\prime}$.

Suppose $a \ge c$.    Let $L$ be the line bundle $V(I_3/I_1I_2)$ 
(respectively $V(I_1/(x^2,y))$.  Then there is a map

$$V(I_2/I_2^2)^c\otimes V(I_1/I_2)^b
\otimes L^{a-c}
\rightarrow V(I^{\prime}/J^{\prime}).$$
By Proposition 2.1 of \cite{R2} it is an isomorphism.  

Similarly, if $c\ge a$, letting $L$ be the line bundle 
$V(I_2/I_3)$ (respectively $V((x^2,y)/I_1^2)$) 
the map 
$$V(I_2/I_1^2)^a\otimes V(I_1/I_2)^d
\otimes L^{c-a}
\rightarrow V(I^{\prime}/J^{\prime})$$
is an isomorphism.  

Thus if the lemma holds for the four line 
bundles that we called $L$ it always holds. The bundles 
$V(I_3/I_1I_2)$ and $V((x^2,y)/I_1^2)$ are the 
tautological bundles over the spaces 
$C_3$ and $C_{2,2}$ respectively and hence have first Chern classes 
$-h_3$ and $-h_3^{\prime}$ respectively.  Thus applying the Whitney 
product formula together with the knowledge of the middle elements from 
Lemma~\ref{basering} to the sequences 

$$0 \rightarrow V(I_3/I_1I_2) \rightarrow V(I_2/I_1I_2) 
\rightarrow  V(I_2/I_3)\rightarrow 0$$
and         
$$0\rightarrow  V((x^2,y)/I_1^2)\rightarrow 
V(I_1/I_1^2)\rightarrow  V(I_1/(x^2,y))\rightarrow 0$$
we see that these four line bundles also satisfy the lemma. $\square$   

\section{Some Chow rings}

In this section we give the Chow rings of some spaces that can be
expressed as fiber bundles over the projectivized cotangent bundle of
$S$ with toric varieties as fibers.  We freely use the material in
\cite{F1}.  Given a ray $r$, we let $v(r)$ denote the smallest
integral point which $r$ passes through.      

\begin{lemma}\label{slick}
Given an exact sequence of vector bundles
$$0\rightarrow V_2 \rightarrow V_1 \rightarrow V_3\rightarrow 0,$$
the projectivization of $V_2$ inside of the projectivization of $V_1$ 
has class $$c_m({\mathcal O}_{V_1}(1) \otimes V_3)$$
where $m$ is the rank of $V_3$. 
\end{lemma}

\noindent{\bf Proof}: The map from $V_1$ to $V_3$ gives a section of 
$\textup {Hom}({\mathcal O}_{V_1}(-1), V_3)$ which vanishes 
on the image of $V_2$.
$\square$

\vspace{.1 in}
The above lemma was communicated to me by Mike Roth.
 
\begin{lemma}\label{intersection}
Let $Y$ be a toric variety of dimension $2$ with corresponding fan
$\Delta$.  Let $D_1$ be a torus
invariant divisor corresponding to a ray $r_1$ in $\Delta$.  Let $D_0$
and $D_2$ be the torus invariant divisors intersecting $D_1$
corresponding to rays $r_0$ and $r_2$ just counter-clockwise and just
clockwise of $r$ respectively.  Then the intersection multiplicity of
$D_0$ and $D_1$ is 
$$D_0D_1= \frac{1}{v(r_1)\wedge v(r_2)}$$
and the self-intersection number of $D_1$ is given by 
$$D_1^2= \frac{v(r_0) \wedge v(r_2)}
{(v(r_0)\wedge v(r_1))(v(r_1) \wedge v(r_2))}.$$
\end{lemma}

\noindent{\bf Proof:}  In the Chow ring of $Y$, one has the relation 
$$\sum_{r_i\in \Delta} (v \cdot v_i)D_i=0$$
for any vector $v$.  Taking $v$ to be orthogonal to $v_0$  
and intersecting with $D_1$ we get the relations 
$$(v_0 \wedge v_2)D_0D_1+(v_1\wedge v_2)D_1^2=0.$$  
Similarly, taking $v$ instead to be orthogonal to $v_2$, we get the
relation
$$(v_2 \wedge v_0)D_2D_1+(v_1\wedge v_0)D_1^2=0.$$ 
It is enough to show that $D_0D_1$ is as claimed because substituting
this into the latter relation, we can solve for the self-intersection
number of $D_1$.  Since $D_2$ does not effect the intersection of
$D_9$ and $D_1$, we can assume that $D_1$ and $D_2$ intersect in a
smooth point.  This implies $D_1D_2= v_2\wedge v_1 = 1$.  Using this
to simplify the relations above, we see that $D_0D_1$ is as claimed.  

\vspace{.1 in}
 
Before stating the main result of this section, we recall some
definitions from \cite{R2}.  Given a sequence of ideals $I_1, \dots ,
I_r$ with 
measuring sequence $m(4,1)$ or $m(3,2)$, recall that  
$U(I_1, \dots, I_r)$ is naturally a fiber over the projectivized
tangent bundle of $S$ with fiber isomorphic to 
$G((x,y^2)/G((x,y^4))$ and $G((x,y^2)/G((x,y^3),(x^2,y))$
respectively.  The normalization of the closure of this fiber is a
toric variety.  We say that the standard fan of this toric variety is
the fan with a ray through $(-1,0)$ corresponding to the divisor
corresponding to automorphisms sending $x$ to $x+ty^2$ for some $t\in
K$ and fixing $y$ and a ray through $(0,-1)$ corresponding to
automorphisms sending $x$ to $x+ty^3$  and fixing $y$ if the measuring
sequence is $m(4,1)$ and fixing $x$ and sending $y$ to $y+tx$ if the
measuring sequence is $m(3,2)$.  Moreover, we say that a ray in a
standard fan is a bounding ray if it corresponds to a boundary divisor
with only one $G((x,y^2))$ fixed point.  It was proved in \cite{R2}
that in all characteristics but $2$, the rays through $(0,1)$ and
$(1,2)$ are bounding rays if they occur and that there are no other
bounding rays.  In characteristic $2$, $(0,1)$ is again a bounding ray
if it occurs, but any other bounding ray must lie in the interior of
the convex cone bounded by the rays through $(1,2)$ and $(0,-1)$.  

\begin{theorem}\label{chowring}
Let $Y$ be an $\aut(R)$ equivariant  
compactification of $U((x,y^4))$ (respectively $U((x,y^3), (x,y^2))$) 
over the space $B = U((x,y^2))$.  Let $\Delta$ be the standard 
fan (as defined in \cite{R2}) 
corresponding to the fiber of $Y$ over $B$.  Label the 
rays in $\Delta$ clockwise starting from the ray through $(-1,0)$ 
so that ${i}^{\rm th}$ ray is labelled $r_{i-1}$.  Let $(n_i, m_i)$ 
be the point of smallest positive distance from 
the origin in $r_i$ having integral coordinates.  
Let $r+2$ be the number of rays in $\Delta$.  For $1\le i \le r$, let 
$D_i$ be the boundary divisor corresponding to $r_i$.  If either the
characteristic of $K$ is not $2$ or there is no bounding ray in the
interior of the cone bounded by the rays through $(1,2)$ and $(0,-1)$
then the Chow ring 
of $A(Y)$ is generated over $A(B)$ by 
the classes of the boundary divisors which we will also denote $D_i$ 
by abuse of notation and the relations are generated by 
$$D_k^2 =s_kD_kD_{k-1} +D_k(a_{k+1}h_2 + b_{k+1}(c_1+h_2))$$ 
for $1 < k \le r,$ 
$$D_k^2 =s_kD_kD_{k+1} - D_k(a_{k-1}h_2 + b_{k-1}(c_1+h_2))$$
for $1 \le k < r,$ and
$$D_iD_j=0$$ for $|i-j|\ge 2$
where $s_k$ is the self-intersection number of the fiber of $D_k$   
and $$(a_k, b_k) =(m_k-n_k,2m_k-3n_k)$$ 
(respectively 
$$(a_k,b_k) = (m_k+n_k,2m_k+n_k)).$$

\end{theorem}

\noindent{\bf Proof}:  The Chow ring of $A$ is generated over the 
Chow ring of $U((x,y^2))$ by the boundary divisors because the restrictions
of the boundary divisors to the fibers generate the Chow rings of the
fibers.  If $|i-j|\ge 2$ and $\{i,j\} \ne \{0,r+1\}$ then the relation  
$D_iD_j=0$ follows from the fact that $D_i$ does not intersect 
$D_j$.  To verify the remaining relations, first we will show that   
$D_i$ is isomorphic to the projectivization of any 
bundle $V(I/J)$ such that $I$ and $J$ have measuring sequence 
at most $(x,y^2), (x,y)$ and $I$ is generated over $J$ by 
$x^{a_i}$ and $y^{b_i}$.  Then we show that 
if $i<r$ then the intersection of $D_{i+1}$ with 
$D_i$ is given by the projectivization of the vector bundle 
$V(I_1/J)$ where $I_1$ is generated over 
$J$ by $x^{a_i}$.  Moreover, we show that 
if $i>1$, then the intersection $D_{i-1}$ 
with $D_i$ is given by the projectivization of the vector  
bundle $V(I_2/J)$  where $I_2$ is generated over 
$J$ by $y^{b_i}$.  Then an application of 
Lemma~\ref{shift} and 
Lemma~\ref{slick} give us some of the relations.  The remaining relations
will be verified through those $Y$ that are projectivizations of staircase 
bundles.            
 
The function $f:\A^2 -0 \rightarrow \P^{1}$
given by $$f((a,b))= (a^{m_i},b^{n_i})$$ extends to a regular function on 
the fiber of $D_i$ which we will also call  $f$.  The map 
$$\varphi: D_i \rightarrow \P V(I/J)$$ such that for a point 
$p$ in the fiber of $D_i$ over $(x,y^2)$, if $f(p) = (s,t)$ then 
$$\varphi (p)= (s x^{a_i} + ty^{b_i})+J$$
is well defined because it is 
independent of the choice of $x$ and $y$.  Moreover, it can be shown to 
be an isomorphism.  If $i<r$, 
the intersection of $D_i$ with $D_{i+1}$ restricted to 
a fiber is $f^{-1}(1,0)$ where $f$ is restricted to the fiber of $D_i$ over 
$(x,y^2)$.  Hence, if $i >1$  
$$\P V(I_1/J)= D_{i+1}D_i $$ 
and if $i<r$  $$\P V(I_2/J)= D_{i-1}D_i.$$  

Let $\xi$ be the hyperplane class of 
$\P V(I/J)$.  By Proposition~\ref{shift} 
$$A(D_i)= A(B)[\xi]/(\xi -a_ih_2)(\xi + b_i(c_1+h_2)).$$   
Lemma~\ref{slick} gives us the relations 
$$(\xi -a_ih_2)D_i = D_iD_{i-1}$$ for $i>1$ and 
$$(\xi +b_i(c_1+h_2))D_i=D_iD_{i+1}$$  
for $i<r$.  Although, there are no global divisors $D_0$ or $D_{r+1}$, if 
$r_i$ is not a bounding ray, there are global divisors $D_iD_{i \pm 1}$ 
Since the intersection of the fiber of $D_i$ with the coordinate axes is 
then $G((x,y^2))$ invariant.  So, we extend these equations to any $i$ 
such that $r_i$ is not a bounding ray.  
Taking the difference of the two equations 
we obtain 
$$D_i(D_{i+1}-D_{i-1}) = D_i(a_ih_2+b_i(c_1+h_2)).$$
  
Hence for $1<i\le r$ and $r_{i-1}$ not a bounding ray  
$$D_{i-1}D_i^2 = D_{i-1}D_i(a_{i-1}h_2+b_{i-1}(c_1+h_2))$$
and for $1\le i <r$ and $r_{i+1}$ not a bounding ray  
$$D_{i+1}D_i^2 = -D_{i+1}D_i(a_{i+1}h_2+b_{i+1}(c_1+h_2)).$$

Since similar relations hold in the Chow ring of the fiber of $Y$ 
over $B$, the Chow ring of $Y$ has relations of the 
form $$D_i^2 = s_iD_iD_{i-1} + D_i\eta_1$$ for $i>1$ 
and $$D_i^2 = s_iD_iD_{i+1} + D_i\eta_2$$  for $i<r$ where 
$s_i$ is the self-intersection number of the fiber $D_i$ over $B$ and 
$\eta$ is the pullback of a class from $B$.  Multiplying these two equations 
by $D_{i+1}$ and $D_{i-1}$ respectively, we see that unless $r_i$ is a 
bounding ray, if $i>1$, 
$$\eta_1 =a_{i+1}h_2+b_{i+1}(c_1+h_2)$$
and if $i<r$,
$$\eta_2 =-a_{i-1}h_2-b_{i-1}(c_1+h_2).$$

\begin{table}\label{chowdata}
\begin{tabular}{|p{.25 in}|p{.75 in}|p{.83 in}|p{.8 in}|p{.57 in}|p{.8
      in}|} 
\hline 
base&rays&$J_0$&$J_1$&$J_2$&$J_3$\\ \hline
$C_3$&$(1,n),(0,1)$&
$\underbrace{2,\dots,2}_{n},3,0$&
$\underbrace{2,\dots,2}_{n+1},1$&
$\underbrace{2,\dots,2}_{n+2}$&
$\underbrace{2,\dots,2}_{n},3,1$ \\ \hline
$C_{2,2}$&$(n,1),(1,0)$&
$\underbrace{1,\dots,1}_{n-1},0,3$&
$\underbrace{1,\dots,1}_{n-1},0,2$&
$\underbrace{1,\dots,1}_{n+1}$&
$\underbrace{1,\dots,1}_{n},2$ \\ \hline
$C_3$&$(1,2),(0,1)$&$4$&$3$&$1,2$&$1,3$\\ \hline
$C_3$&$(1,n),(0,1)$&
$\underbrace{2,\dots,2}_{n-3},1,4$&
$\underbrace{2,\dots,2}_{n-3},1,3$&
$\underbrace{2,\dots,2}_{n-1}$&
$\underbrace{2,\dots,2}_{n-2},3$\\ \hline

\end{tabular}
\caption{Table of Projective Bundles}
\end{table}

It remains to verify the relations involving $D_i^2$ for $r_i$ a bounding 
ray.  Since the relations depend on the neighborhood of $D_i$, it is
enough 
to verify them for some space $Y$ for each possible pair of  
rays corresponding to adjacent boundary divisors, 
such that one of the rays is a bounding ray.  Any such 
pair occurs for a space $Y(i,n)$ given by the $i^{\rm th}$ row of 
Table 1 as follows.  The integer $n$ must be at least $i-1$, 
except that 
$Y((n,3))$ is independent of $n$.  
Let the ideal $J_k$ be as given by the 
entries in the $i^{\rm th}$ row of the respective column.  We define 
the space $Y(n,i)$ to be the projectivization of the staircase 
bundle $V(J_1/J_3)$ over the base $B^{\prime}$ given in the first column.  By 
Theorem~\ref{proj}, this is also the space obtained by superimposing 
$B^{\prime}$ with the space $C(J_0)$.  
Thus if $i=3,4$ (respectively $i=1,2$) then 
$Y(n,i)$ is a compactification of     
$U((x,y^4))$ (respectively $U((x,y^3), (x^2,y))$.
The second entry in the $i^{\rm th}$ row gives the rays in the fan of the 
cone of the fiber of $Y(n,i)$ over $B$ corresponding to the two boundary 
divisors.  The first ray corresponds to the pullback of 
the boundary of $B^{\prime}$ and the second to the projectivization of the 
sub-bundle $V(J_2/J_3)$ over $B^{\prime}$.  By Lemma~\ref{slick}, these 
two divisors have classes $h_4 + c_1(V(J_1/J_2))$ and $h_3-h_2$ (respectively 
$h_3+h_2+c_1$) where $h_4$ is the hyperplane class of $Y(n,i)$.  
Using this to eliminate $h_3$ and $h_4$ in 
the relation from the base $C_3$ or $C_{2,2}$ 
and the relation 
$$(h_4+c_1(J_1/J_2))(h_4+c_1(J_2/J_3))= 0$$
we recover the relations 
in the statement of the theorem.  $\square$  

\begin{proposition}\label{classes} Keeping the notation of the proof of  
Theorem~\ref{chowring}, let $J_0$, $J_1$, $J_2$, and $J_3$ be the
ideals associated to the space $Y(i,n)$ as given in
Figure~\ref{chowdata} and $r_1$ and $r_2$ the two rays given in the
$i^{\rm th}$ column of the table.  Let $D$ be the
divisor of $Y(i,n)$ corresponding to the ray through $(1,0)$ if $i=2$ 
and the ray through $(0,1)$ otherwise.  Then 
$$c_1(V(J_0/J_3))= -D+ c_1(J_1/J_2)$$ and 
$$c_1(V(J_1/J_0))= D+ c_1(J_2/J_3).$$   
\end{proposition}

\noindent{ \bf Proof:}  The bundle $V(J_0/J_3)$ is the tautological
bundle over $Y(i,n)$ and thus has first Chern class $-h_4$.  Recalling
from the proof of Theorem~\ref{chowring} that $$D(r_2) =
h_4+c_1(V(J_1/J_2))$$ 
and the fact that 
$$ c_1(J_1/J_0)+ c_1(J_0/J_3) = c_1(J_1/J_2)+ c_1(J_2/J_3),$$
the proposition follows.  $\square$   
 
\section{Examples}
In this section we find the number $N_i$ of curves  
in a suitable linear series $\L$ on a surface $S$ with the singularity with
local equation $x^2 +y^i$ for $i$ from $2$ to $8$.  We will let $D$
denote the divisor of a section of $\L$ and $L$ the associated line
bundle on $S$.  Moreover, we will let $I_m$ denote the ideal
$(x,y^m)$.  The ideal 
$$B_i= (x^2, xy^{\left [ \frac{i}{2} \right ]}, y^i)$$
corresponds to this type of singularity in the sense of section 3.
The linear series $\L$ will be of projective dimension $n-1$, giving 
$n$ independent sections of the bundle $V_L(B_i)$ over the space 
$C(B_i)$.   Up to scaling, the linear combinations of these sections 
with zeroes are in bijection with curves in $L$ with the given
singularity.  Thus the number of these curves is the number of the 
vector bundle $V_L(B_i)$.  This Chern number can be found by 
finding the Chern class of dimension $0$ of $V_L(B_i)$ and then using 
the Chow ring of the space $C(B_i)$ or some other space dominating 
to find the degree of this class.  Up to $i=6$, this information is
given by Lemma~\ref{basering} and Proposition~\ref{shift}.  For 
$i=7$ and $i=8$, the relevant Chow rings are given by
Theorem~\ref{chowring}, but we will have to use slightly add-hoc means
to find the Chern classes we are interested in.    
 
To find the relevant Chern classes, we will find the Chern classes of
the successive quotients in the sequence 
$$V(I_4^2) \rightarrow V(I_4^2+I_1^7)\rightarrow 
V(I_3^2)\rightarrow V(I_3^2+I_1^5) \rightarrow 
\break V(I_2I_3) \rightarrow $$
$$V(I_2^2) 
\rightarrow
V(I_2^2+I_1^3)\rightarrow \break V(I_1I_2)\rightarrow
V(I_1^2)\rightarrow V(I_1)$$  
and then apply the Whitney product formula together with
Lemma~\ref{tensor}.

To find the number $N_2$, the number of nodes in a pencil of curves,
we need only work over the surface $S$.  Thus $N_2$ is equal to the 
second Chern class of $V(I_1^2)\otimes L$.  Since this Chern class is 
is already expressed in terms of the Chern classes of the surface and 
and the divisor $D$, no further substitution is necessary.  This
example as well as the example of computing $N_3$ are worked out in
detail in \cite{V}.  To find 
$N_3$ and $N_4$, the numbers of cuspidal and tacnodal curves in a
suitable linear series $\L$, we will work over the projectivized
cotangent bundle of $S$.  The Chern classes $c_3(V(B_3)\otimes L)$ and 
$c_3(V(B_4) \otimes L)$ come expressed in terms of divisors pulled back from
$S$ and the hyper-plane class $h_2$ of the cotangent bundle.  Using
the relation $h_2^2+ c_1h_2+c_2=0$ to make these classes linear in
$h_2$, the numbers, $N_3$ and $N_4$ are the coefficients of $h_2$ of
these classes.  Similarly, to find $N_5$ and $N_6$, we work over the
space $Y(I_2, I_3)$.  The fourth Chern classes of the bundles 
$V(B_5) \otimes L$ and $V(B_6) \otimes L$ come expressed in terms of
pullbacks of divisors on $S$, the pullback of the hyperplane class
$h_2$ of the projectivized cotangent bundle of $S$ and the hyperplane 
class $h_3$ of $C(I_1, I_2)$.  Using the relations in the Chow ring of
$C(I_1, I_2)$ to make these Chern classes linear in $h_3$ and $h_2$ 
separately, the numbers $N_5$ and $N_6$ are the coefficients of 
$h_2h_3$ in these two classes.  

Finding the numbers $N_7$ and $N_8$ is
a bit trickier because we must find the 
Chern polynomials $c(V(B_7) \otimes L)$ and
$c(V(B_8)\otimes L)$ by slightly add-hoc means.  In particular, we
need to find the Chern classes of the line bundles 
$V(I_3^2/I_3I_4)$, $V(I_3I_4/B_7)$, and $V(B_7/B_8)$.  
We will use Table~\ref{vanish} to see where the maps 
$$\varphi_1: V(I_3/I_4)^2 \rightarrow V(I_3^2/I_3I_4),$$
$$\varphi_2: V(I_2/I_3)^3 \rightarrow V(I_3^2/I_3I_4),$$
$$\varphi_3: V(I_4/I_1I_3)\otimes V(I_3/I_4) \rightarrow
V(I_3I_4/B_7),$$
and $$\varphi_4: V(I_1I_3/I_1I_4)\otimes V(I_3/I_4) \rightarrow 
V(B_7/B_8)$$
and then apply Porteous's formula.  We will work over the
compactification $Y$ of $U(I_4)$ with boundary divisors corresponding
to rays through points $(0,1)$, $(1,4)$, $(1,3)$, $(2,5)$ and $(1,2)$
in the standard fan of the fiber, since this is the smallest bundle
over which all of the vector bundles we will use are defined.  Some
vector bundles are also defined over smaller spaces.      
The Chern classes over the smaller spaces can be related the 
Chern classes over $Y$ via the following lemma.     

\begin{lemma}\label{pullback}  Let $Z$ and $Z^{\prime}$ be two toric 
varieties of dimension $2$ with fans $\Delta$ and $\Delta ^{\prime}$ 
such that $\Delta ^{\prime}$ is a subdivision of $\Delta $.  
With respect to 
the map from $Z^{\prime}$ to $Z$ compatible with these fans, the
pullback of a divisor $D$ corresponding to a ray $r_1$ in $\Delta $ is
of the form 
$$\sum_i a_i D_i $$ 
where $i$ indexes of the rays $r_i$ in $\Delta ^{\prime}$ and the 
$D_i$ 's are the corresponding divisors.  The integers $a_i$ can be
found as follows.  If $r_i$ lies strictly between $r_1$ and an adjacent
ray $r_2$ in $\Delta$, then if 
$$v(r_i)= l_1v(r_1) +l_2v(r_2)$$ 
then $a_i$ is equal to $l_1$.  If $r_i= r$ then $a_i = 1$.  
In any other case, $a_i =0$.  
\end{lemma}

\noindent{\bf Proof:}  It is enough to check the lemma in the case
that  $\Delta ^{\prime}$ differs $\Delta $ by a single subdivision
corresponding to a ray $r_3$ because the coefficient of the divisor of
corresponding to this ray is the same as the coefficient of this ray
for any $\Delta ^{\prime}$ containing $r_3$.  The self-intersection
number of $D$ is the same as the self-intersection number of its
pullback.  Writing the pullback of $D$ as $D_1 + qD_3$ we have 
$$D^2= D_1^2+ 2qD_1D_3 + q^2D_3^2.$$
Thus if $D_1^2= D^2$, then $q=0$.  This verifies the lemma in the cases
when $r_1= r_3$ and when $D_1$ and $D_3$ do not intersect.  It remains
to verify the lemma when $r_3$ subdivides a cone bounded by 
$r_1$ and another ray $r_2$.  Writing $r_3 = q_1 r_1 + q_2r_2$ and
substituting, using Lemma~\ref{intersection} the above equation
reduces to $(q - q_1)^2 = 0$.  Thus $q= q_1$ and we have verified the
lemma.  $\square$

\vspace{.1 in}

{\tiny
\begin{table}\label{vanish}
\begin{tabular}{p{.11in} |p{.7in}|p{.7in}|p{.9in}|p{.8in}|p{1in}|}
&$(0,1)$&$(1,4)$&$(1,3)$&$(2,5)$&$(1,2)$\\ \hline

$0,3$&$(x+ay^2)+ I_1I_2$&$I_1^2$&$I_1^2$&$I_1^2$&$I_1^2$ \\ \hline

$0,4$&$I_1I_2$&$I_1I_2 $&$I_1I_2 $&$I_1I_2$&$({a^2y^2- bxy})+ 
\break I_1^3+I_2^2$ \\ \hline

$1, 3$&$ ({xy + ay^3}) + I_2^2$&$I_1^3+I_2^2$&$I_1^3+I_2^2$&
$I_1^3+I_2^2$&
$I_1^3+I_2^2$ \\ \hline

$4, 3$& $y(x + ay^2)^2 +I_2^3$&$I_2^3+I_1^5$&
$({bx^3+2a^3xy^3}) +\break I_1^2I_2^2 + I_1^5 $&
$I_1^3I_2$&$({a^2y^4- 2bxy^3})+  \break I_1^2I_2^2 +I_1^5$
\\ \hline

$3, 4$&$ I_2(xy+ay^3)+ I_2^3 $&
$I_1^2I_2$&$I_1^2I_2$&
$({3b^2x^3+a^5y^4})+ \break I_1^3I_2
$& $I_1^4$ \\ \hline


$1, 4$&$I_2^2$&$I_2^2$&$({2bx^2+a^3y^3})+\break I_1^2I_2 $&
$I_1^3$&$I_1^3$ \\ \hline

$4,4$&$I_2^3$&$(bx^3+a^4y^5)+ \break I_1^2I_2^2$&$I_1^2I_2^2+I_1^5$&
$I_1^2I_2^2+I_1^5$&$(x^2,({bx-2ay})^2)^2+ \break I_1^5$ \\ \hline

\end{tabular}
\caption{Boundary ideals}
\end{table}}

\normalsize 

\begin{table}\label{c1}
\begin{tabular}
{p{.5in} |p{.12in}|p{.15in}|p{.4in}|p{.4in}|p{.4in}|p{.4in}|p{.4in}|}
& $c_1$& $h_2$ & $D(0,1)$ & $D(1,4)$ & $D(1,3)$ & $D(2,5)$ & $D(1,2)$
\\ \hline
1,1/2,2&2&2&&&&& \\ \hline
1,2/2,1&1& &&&&& \\ \hline
2,1/2,2&3&3&&&&& \\ \hline
2,2/2,3&4&4&&2&2&4&2 \\ \hline
2,3/3,2&2&1&&&&& \\ \hline
3,2/3,3&5&5&&2&2&4&2 \\ \hline
3,3/3,4&6&6&2&6&4&6&3 \\ \hline
3,4/4,3&3&2&&&&1&\\ \hline
4,3/4,4&7&7&2&6&5&8&3 \\ \hline
0,2/0,3&2&2&&1&1&2&1 \\ \hline
0,3/1,2& &-1&&-1&-1&-2&-1 \\ \hline
1,2/1,3&3&3&&1&1&2&1 \\ \hline
0,3/0,4&3&3&1&3&2&3&1\\ \hline
0,4/1,3& &-1&-1&-3&-2&-3&-1\\ \hline
1,3/1,4&4&4&1&3&2&4&2\\ \hline

\end{tabular}
\caption{First Chern classes}
\end{table}

Each column of Table 2 corresponds to a boundary divisor of
$Y$ and each 
row corresponds to an ideal as labelled on the top row and left-most 
column respectively.  The sequence $m,n$ corresponds to the ideal
$x^2, xy^m, y^{m+n}$.  Let $J(v, I)$ be the entry corresponding to the ray 
through $v$ and ideal $I$.  This ideal is the projection of the image 
of $I$ under the generic automorphism 
to the smallest plane in the Pl{\" u}cker embedding of the fiber $C(B,I)$ 
over $B$  
containing the boundary divisor corresponding to $v$.   

The map $\varphi_1$ drops rank on the boundary divisors 
$D_i$ such that 
$$J(v_i,I_3)^2 \subset J(v,I_3I_4).$$  
Hence by Porteous's formula over $Y$ we have 
$$c_1(V(I_3^2/I_3I_4)) =
2c_1(V(I_3/I_4))+ n_1D(1,2)$$
$$= 
6(c_1+h_2)+2D(0,1)+ 6D(1,4)+4D(1,3)+6D(2,5)+ (2+n_1)D(1,2)$$
where $n_1$ is a positive integer.  
Here we have found $c_1(V(I_3/I_4))$ through Proposition~\ref{classes}
together with Lemma~\ref{pullback}.  Similarly, one can find the Chern
classes of $V(I_4/I_1I_3)$ and $V(I_1I_3/I_1I_4)$ over $Y$ in this
way.  Following the same steps, Porteous's formula applied
to the remaining $\varphi_i$'s yields 
$$c_1(I_3^2/I_3I_4)= 6(c_1+h_2)+ n_2D(0,1)+ (n_3+3)D(1,4) +
(n_4+3)D(1,3)+ $$
$$6D(2,5)+ 3D(1,2),$$  
$$c_1(I_3I_4/B_7)=3c_1+2h_2 + n_5D(2,5),$$
and $$c_1(B_7/B_8)= 7(c_1+h_2)+ 2D(0,1)+ 6D(1,4)+(4+n_6)D(1,3)+ $$
$$(7+n_7)D(2,5)+3D(1,2).$$  
Comparing the two expressions for $c_1(V(I_3^2/I_3I_4)$ we see that $n_1= 1$.
Since $C(I_2, I_3^2, B_7)$ does not have a boundary divisor $D(2,5)$,
by Lemma~\ref{pullback} the coefficient of $D(2,5)$ in the first Chern
class of $V(I_3^2/B_7)$ as a vector bundle over $Y$ is 
the sum of the coefficients of $D(1,3)$ and $D(1,2)$.  By the Whitney
product formula 
$c_1(V(I_3^2/B_7))$ is the sum of $c_1(V(I_3^2/I_3I_4))$ 
and $c_1(V(I_3I_4/B_7))$, we obtain $n_5=1$.  Similarly, since 
$C(I_2, I_3^2,I_4^2)$ does not have boundary divisors 
$D(1,3)$ and $D(2,5)$, the coefficient of $D(1,3)$ in 
$c_1(V(I_3^2/I_4^2)$ as a class in $Y$ 
is half the sum of the coefficient of $D(1,4)$
and the coefficient of $D(1,2)$.  Moreover, the coefficient of
$D(2,5)$
is half the sum of the coefficients of $D(1,4)$ and three times the
coefficient of 
$D(1,2)$.  By the Whitney product formula, we have 
$$c_1(V(I_3^2/I_4^2)= c_1(V(I_3^2/I_3I_4))+ c_1(V(I_3I_4/B_7)+ 
c_1(V(B_7/I_4^2)).$$  Thus it follows that $n_6=n_7=1$.  
The first Chern classes we have found or used along the way are listed
in Table 3.  The leftmost entry $m_1,m_2/ m_3, m_4$ of 
each row signifies that the remaining entries in that row are the
coefficients of the divisors listed at the top of each column 
in the Chern class 
$c_1(V((x^2, xy^m_1, y^{m_1+m_2}))/(x^2, xy^m_3,
y^{m_3+m_4})))$ over $Y$.  

Having found all relevant Chern classes, using the Whitney product
formula together with Lemma~\ref{tensor}, we can find the 
Chern classes $c_5(V(B_7)\otimes L)$ and $c_5(V(B_7)\otimes L)$.  It
remains to calculate the degrees of these Chern classes using the
relations for the Chow ring of $Y$ given in Theorem~\ref{chowring}.
Using these relations, these classes can be made free of squares of
boundary divisors and linear in $h_2$.  The respective degrees of the 
Chern classes are then the sums of the coefficients of terms of the
form $h_2D_iD_{i+1}$.  Our results are summarized below.  

$$N_ 2=6 D^2 + 4 D c_1 + 2 c_2 $$ 
$$N_3 =12 D^2 + 12 D c_1 + 2c_1^2 + 2 c_2$$ 
$$N_4 =50 D^2 + 64 Dc_1 + 17c_1^2 + 5c_2$$ 
$$N_5 =180 D^2 + 280Dc_1 + 100c_1^2$$ 
$$N_6 =630D^2 + 1140Dc_1 + 498c_1^2 - 60c_2$$ 
$$N_7 =2128D^2 + 4368Dc_1 + 2232c_1^2 - 424c_2$$
and 
$$N_8 =7272D^2 + 16544Dc_1 + 9548c_1^2 - 2148c_2.$$

\end{document}